\documentclass{IEEEtran4PSCC}
\usepackage{xcolor}
\usepackage{amsfonts,subcaption,multirow, booktabs}
\usepackage{amsthm}
\newtheorem{remark}{Remark}
\newtheorem{assumption}{Assumption}
\newcommand{\bs}{\boldsymbol}
\newcommand{\mc}[1]{\mathcal{#1}}

\ifCLASSINFOpdf
   \usepackage[pdftex]{graphicx}
\else
   \usepackage[dvips]{graphicx}
\fi
\usepackage[cmex10]{amsmath}
\allowdisplaybreaks
\makeatletter
\let\old@ps@headings\ps@headings
\let\old@ps@IEEEtitlepagestyle\ps@IEEEtitlepagestyle
\def\psccfooter#1{%
    \def\ps@headings{%
        \old@ps@headings%
        \def\@oddfoot{\strut\hfill#1\hfill\strut}%
        \def\@evenfoot{\strut\hfill#1\hfill\strut}%
    }%
    \def\ps@IEEEtitlepagestyle{%
        \old@ps@IEEEtitlepagestyle%
        \def\@oddfoot{\strut\hfill#1\hfill\strut}%
        \def\@evenfoot{\strut\hfill#1\hfill\strut}%
    }%
    \ps@headings%
}
\makeatother

\psccfooter{%
        \parbox{\textwidth}{\hrulefill \\ \small{23rd Power Systems Computation Conference} \hfill \begin{minipage}{0.2\textwidth}\centering \vspace*{4pt} \includegraphics[scale=0.06]{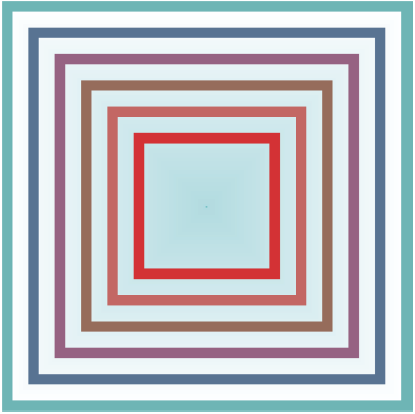}\\\small{PSCC 2024} \end{minipage} \hfill \small{Paris, France --- June 4 -- June 7, 2024}}%
}

\begin{document}
%
\title{Linearized Optimal Power Flow for {\color{black}Multiphase Radial Networks} with Delta Connections}
 
\author{
\IEEEauthorblockN{Geunyeong Byeon and Minseok Ryu}
\IEEEauthorblockA{School of Computing and Augmented Intelligence\\
Arizona State University\\
Tempe, AZ, USA \\
geunyeong.byeon@asu.edu, minseok.ryu@asu.edu}
\and
\IEEEauthorblockN{Kibaek Kim}
\IEEEauthorblockA{Mathematics and Computer Science Division\\
Argonne National Laboratory\\
Lemont, IL, USA\\
kimk@anl.gov}
}


\maketitle

\begin{abstract}
This paper proposes a linear approximation of the alternating current optimal power flow problem for multiphase distribution networks with voltage-dependent loads connected in both wye and delta configurations. We establish a set of linear equations that exactly describes the relationship between power changes at a bus and the corresponding values from a delta-connected device under specific assumptions necessary for a widely accepted linear model. Numerical studies on IEEE test feeders demonstrate that the proposed linear model provides solutions with reasonable error bounds efficiently, as compared with an exact nonconvex formulation and a convex conic relaxation. Our experiments reveal that modeling delta-connected, voltage-dependent loads as if they are wye-connected can lead to significantly different voltage profiles. We also investigate the limitations of the proposed linear approximation stemming from the underlying assumptions, while emphasizing its robust performance in practical situations.
\end{abstract}


\thanksto{\noindent This material is based upon work supported by the U.S. Department of Energy, Office of Science, under contract number DE-AC02-06CH11357.}

\section{Introduction}
\IEEEPARstart{O}{ptimal} power flow (OPF) refers to an optimization problem that aims to find an optimal operating point for a power grid while abiding by physical laws governing the electric power flow as well as various operational limits imposed on network components. OPF is often embedded in a variety of optimization problems that arise in the planning and operation of power grids (e.g., unit commitment and expansion planning problems), but the power flow equations that formulate the steady-state physics of power flow make OPF both nonconvex and NP-hard even for radial networks \cite{lehmann2015ac}. Because of OPF's computational complexity and wide applicability, the past decade has seen an active line of research that attempts to reformulate OPF into a computationally tractable problem, especially as a convex relaxation and/or a linear approximation. For a comprehensive review of various OPF formulations, we refer readers to \cite{molzahn2019survey}.

Most of the literature has focused on transmission grids that feature balanced OPF that can be treated like a single-phase network. OPF for distribution grids has received less attention because they are usually characterized by unidirectional power flow over a huge radial network, which makes their operations largely based on monitoring and controlling a small subset of network components, rather than by solving OPF. However, OPF in  distribution grids becomes important as an increasing number of distributed energy resources introduce bidirectional power flow, which requires active controls based on OPF. Moreover, the single-phase models may not be directly applicable because of the 
unbalanced nature 
of distribution grids, which calls for tractable multiphase OPF formulations for distribution grids. In addition, it is important to model the voltage dependency of loads in distribution grids in order to employ conservation voltage reduction for reducing consumption \cite{schneider2010evaluation}. Given the increasingly important role of OPF in distribution grids,  ongoing efforts seek to develop an open-source framework that facilitates the implementation of OPF for distribution grids \cite{rigoni2020open,fobes2020powermodelsdistribution}. 

For the unbalanced, multiphase OPF, semidefinite programming (SDP) relaxations  have been proposed for networks with  wye-connected constant-power loads  \cite{dall2013distributed,gan2014convex}. A linear approximation has also been proposed by ignoring line losses and assuming nearly balanced voltages  \cite{gan2014convex}. The SDP relaxation is extended in \cite{zhao2017optimal} to
networks with both wye- and delta-connected devices by introducing an additional positive semidefinite matrix that represents the outer product of bus voltages and phase-to-phase currents in the delta connection. In \cite{zhou2021exactness} the authors proposed an algorithm for achieving exactness of this SDP relaxation that adds a penalty term in the objective function. For voltage-dependent loads, approximate representations of ZIP loads were proposed in \cite{usman2020bus} and \cite{molzahn2014approximate}. In this paper we focus on multiphase OPF for networks with delta-connected, exponential loads. The most related work is \cite{claeys2021voltage}, which proposed a convex conic relaxation of exponential loads using power cones that can be added to an SDP relaxation of the multiphase OPF problem.

Convex conic relaxations have several advantages in that they provide lower bounds or even a globally optimal solution (for constant load models) to the nonconvex OPF problems and also are proven to be efficiently solvable by interior-point methods. In practice, however, they often suffer from numerical issues and show slow computation time when solving large-scale problems. This makes the convex conic relaxations difficult to be embedded in planning and operations problems with more features, for example, line switching and transformer controls. 

Alternatively, we propose a linear approximation of multiphase OPF with exponential load models configured as either wye or delta connections.
The main contribution of this paper lies in establishing a set of linear equations that exactly describes the relationship between the power withdrawal/injection from a bus and that from a delta-connected device under the assumptions required by the linear model of \cite{gan2014convex}. 
We demonstrate numerically that the proposed model yields solutions with reasonable error margins (e.g., less than 1\% deviation in voltage magnitude) and shorter computation times compared to the convex relaxation, across three IEEE distribution feeders. {\color{black}Furthermore, a numerical experiment underscores the significance of the proposed model by showing that a nonlinear model neglecting delta-connected, voltage-dependent loads may yield voltage profiles notably different from those of alternating current (AC) OPF compared to the proposed linear model. Additionally, it underscores the robust performance of our proposed model even when the assumptions of low line losses and balanced voltage are reasonably violated, while identifying potential limitations of the proposed model for systems with very high voltage unbalance. }

This paper is organized as follows. Section \ref{sec:prelim} defines preliminaries on the network model, formalizes the AC OPF problem with only wye-connected, constant-power loads, and describes its linear approximation. In Section \ref{sec:vol-dependent-load-model}, their extensions to the networks with delta-connected, voltage-dependent loads are explained and proposed. Section \ref{sec:experiments} reports numerical results, and Section \ref{sec:conclusion} concludes the paper.

\section{Preliminaries}\label{sec:prelim}
\begin{table}[!t]\fontsize{9}{9}\selectfont
\renewcommand{\arraystretch}{1.2}
	\centering
	\caption{Parameters} \label{table:param}
	{\begin{tabular}{p{0.26\linewidth} p{0.66\linewidth}}
		\toprule
		Notation & Description\\
		\midrule
        $(\mc N , \mc E)$ & an undirected graph representing a distribution grid\\
        \ $\mc N$ & set of buses, indexed with $\{0, 1, \cdots, n\}$, where $0$ denotes the substation bus that serves as the slack bus\\
        \ $\mc E$ & set of lines, indexed by a tuple $(e, i, j)$ with $i < j$, where $e$ denotes an integer for indexing the line and $i, j \in \mc N$ denote its two end buses\\
        $\mc N_{{leaf}} \subseteq \mc N$ & set of leaf buses\\
        $\mc P = \{a, b, c\}$ &  three phases of the network\\
        $\mc P_x \subseteq \mc P$ & phases of a network component $x$, where $x$ can be an index of a network component\\
        $\mc B, \mc G,$ $\mc L$ & set of all shunt devices, generators, and loads, respectively\\
        $\mc B(i)$, $\mc G(i)$, $\mc L(i)$ & set of all shunt devices, generators, and loads connected at $i \in \mc N$\\
        $\bs V_0^{ref}$ & reference complex voltage specified for the slack bus $0 \in \mc N$\\ 
        $\underline{v}_{i\phi}$ (and $\overline{v}_{i\phi}$) & lower (and upper) bound of voltage magnitude squared at $i \in \mc N$ on $\phi \in \mc P$\\
        \multicolumn{2}{l}{For each line $(e,i,j) \in \mc E$:}\\
        \ $\bs Y^{sh}_{eij}$ (and $\bs Y^{sh}_{eji})$ & shunt admittance matrices of $e$ near $i$ (and $j$)\\ 
        \ $\bs Z_e$ & series impedance matrix of $e$\\
        \multicolumn{2}{l}{For each shunt device $s \in \mc S$:}\\
        \ $\bs Y^{sh}_s$ & admittance matrix of $s$\\
        \multicolumn{2}{l}{For each generator $k \in \mc G$:}\\
        \ $\mc S_k$ & feasible region of power output of $k$; e.g., $\{\bs p_k + \mathrm{i} \bs q_k:  \underline{\bs p}_k \le \bs p_k \le \overline{\bs  p}_k, \underline{\bs q}_k \le \bs q_k \le {\color{black}\overline{\bs q}_k}\}$ for some $\underline{\bs p}_k$, $\overline{\bs p}_k$, $\underline{\bs q}_k$, $\overline{\bs q}_k$$\in \mathbb R^3) $. \\
        \multicolumn{2}{l}{For each load $l \in \mc G$:}\\
        \ $\bs S_l^0$ & nominal complex power load of $l$ \\
        \bottomrule
    \end{tabular}}
\end{table}

\begin{table}[!t]\fontsize{9}{9}\selectfont
\renewcommand{\arraystretch}{1.2}
	\centering
	\caption{Variables} \label{table:variables}
	{\begin{tabular}{p{0.29\linewidth} p{0.62\linewidth}}
		\toprule
		Notation & Description\\
		\midrule
        \multicolumn{2}{l}{For each bus $i \in \mc N$:}\\ 
        \ $\bs V_{i} = [ V_{i\phi} ]_{\phi \in \mc P}$ & bus voltage at $i$\\ 
        \ $\bs I^{sh}_s = [I^{sh}_{s\phi}]_{\phi \in \mc P}$ & current injection to $i$ from $s \in \mc B(i)$\\
        \ $\bs I^{g}_k = [I^g_{k\phi}]_{\phi \in \mc P}$  & current injection to $i$ from $k \in \mc G(i)$\\
        \ $\bs I^{b}_l = [I^ b_{l\phi}]_{\phi \in \mc P}$ & current withdrawal from $i$ to $l \in \mc L(i)$\\
        \ $\bs S^g_k = [S^g_{k\phi}]_{\phi \in \mc P}$ & power injection at $i$ from $g \in \mc G(i)$\\
        \ $\bs S^b_l = [S^b_{l\phi}]_{\phi \in \mc P}$ & power withdrawal from $i$ to $l \in \mc L(i)$\\
        \multicolumn{2}{l}{For each line $(e,i,j) \in \mc E$:}\\
        \ $\bs I_{eij} = [I_{eij\phi}]_{\phi \in \mc P}$ & current flowing through $e$ from $i$\\
        \ $\bs I^s_{eij} = [I^s_{eij\phi}]_{\phi \in \mc P}$ & series current flowing through $e$ from $i$\\
        \ $\bs I^{sh}_{eij} = [I^{sh}_{eij\phi}]_{\phi \in \mc P}$ & shunt current flowing through $e$ from $i$\\
        \ $\bs S_{eij} = [S_{eij\phi}]_{\phi \in \mc P}$ & series power flowing through $e$ from $i$\\
        \ $\bs I_{eji}$,$\bs I^s_{eji}$,$\bs I^{sh}_{eji}$,$\bs S_{eji}$ & those from its end-bus $j$\\
        \multicolumn{2}{l}{For each load $l \in \mc L$:}\\
        \ $S^d_{l\phi}= p^d_{l\phi} + \mathbf{i} q^d_{l\phi}$ & power consumption of $l$ on  $\phi \in \mc P$\\
        \ $V_{l\phi}$ & voltage applied to $l$ on $\phi \in \mc P$\\
        \ $I_{l\phi}^d$ & current passing through $l$ on $\phi \in \mc P$\\
        \bottomrule
    \end{tabular}}
\end{table}

\begin{figure}
  \centering
  \includegraphics[width=0.85\linewidth]{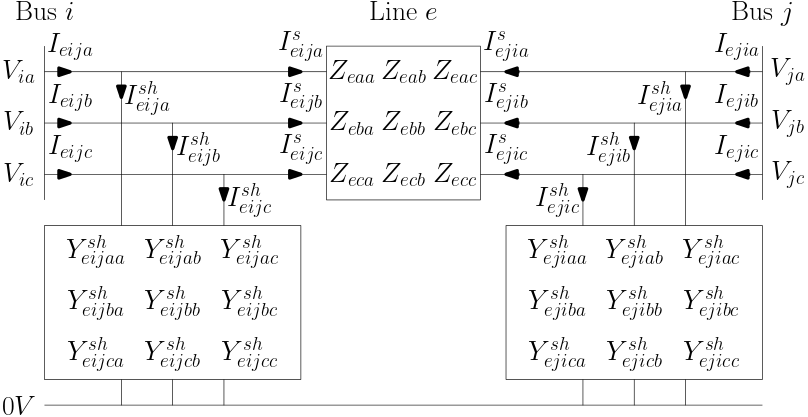}
  \caption{Notations for line $(e,i,j) \in \mc E$}
  \label{fig:line-notation}
\end{figure}

\begin{figure}
  \centering
  \includegraphics[width=0.4\linewidth]{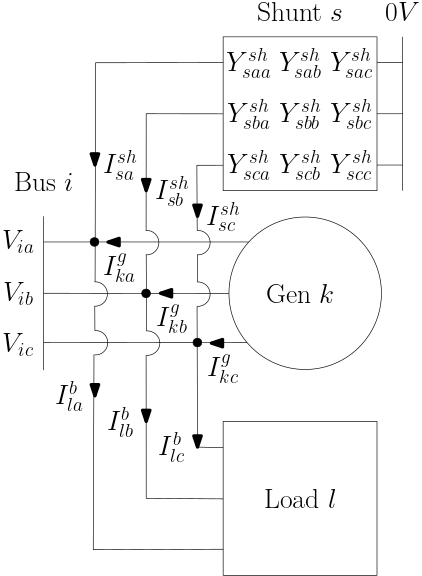}
  \caption{Notations for bus $i \in \mc N$}
  \label{fig:bus-notation}
\end{figure}

The parameters and variables are summarized in Tables \ref{table:param}, \ref{table:variables} and Figures \ref{fig:line-notation}, \ref{fig:bus-notation}, respectively. Throughout this paper, calligraphic letters denote sets, and bold-faced letters denote vectors or matrices depending on the context. Complex numbers are denoted by non-bold, upper-case letters. For a matrix $\bs A$, $\mbox{diag}(\bs A)$  denotes a vector constructed by the diagonal elements of $\bs A$; we abuse notation and use $\mbox{diag}(\bs a)$ for a vector $\bs a$ to represent a matrix in which its diagonal is composed of $\bs a$ and off-diagonal elements are all zeros. For a complex vector or matrix $\bs A$, $\bs A^H$ denotes its conjugate transpose.

Most distribution grids involve three
phases and are operated in radial (i.e., tree) structures, so we
focus on the case where the graph is a tree (i.e., ($\mc N, \mc E$) is connected and $|E| = n$). We consider each phase $a,b,c$ an integer $0,1,2$, respectively, so that for $\phi \in \mc P$, 
(i) $\phi^+$ refers to $(\phi + 1) \bmod 3$ and (ii) $\phi^- := (\phi + 2) \bmod 3$. For example $a^+ = b$ and $a^- = c$. The electric power $\bs S^b_l$ withdrawn from bus $i$ by load $l \in \mc L(i)$ is a function of the bus voltage, say, $f_l(\bs V_{i})$. The form of the function $f_l$ depends on how the load is modeled (e.g., constant power, constant impedance, exponential load) and configured (e.g., delta or wye). In this preliminary section, $f_l(\bs V_i)$ is considered as a constant function $f_l(\bs V_{i}) = {\color{black}\bs S_l^0}$ for its nominal power ${\color{black}\bs S_l^0}$. \emph{In Section \ref{sec:vol-dependent-load-model}, we consider a more general form of $f_l(\bs V_{i})$ that can model delta-connected voltage-dependent load.}

Power flows are governed by the following physical laws:
\begin{itemize}
  \item[(i)] Due to Ohm's law, 
 $\bs V_i - \bs V_j = {\color{black}\bs Z_{e}} \bs I^s_{eij}$, $\bs I^{sh}_{eij} = \bs Y^{sh}_{eij} \bs V_i$, $\bs I^{sh}_{eji}=\bs Y^{sh}_{eji} \bs V_j$ for $(e,i,j) \in \mc E$ and ${\color{black}\bs I^{sh}_{s}=-\bs Y^{sh}_{s}} \bs V_i$ for $i \in \mc N$ and $s \in \mc B(i)$. 
  \item[(ii)] Due to Kirchhoff's current law, for $i \in \mc N$, we have $\sum_{(e,i,j) \mbox{ or } (e,j,i) \in \mc E} (\bs I^s_{eij} + \bs I^{sh}_{eij}) =  \sum_{s \in \mc B(i)} \bs I^{sh}_{s} + \sum_{k \in \mc G(i)} \bs I^{g}_{k} - \sum_{l \in \mc L(i)} \bs I^{b}_l$. For $(e,i,j) \in \mc E$, $\bs I^s_{eij} = - \bs I^s_{eji}$.
  \item[(iii)] Due to the power equation, 
  for $(e,i,j) \in \mc E$ and $\phi \in \mc P$, $\bs S_{eij\phi} = V_{i\phi} (I^s_{eij\phi})^H$ and $\bs S_{eji\phi} =V_{j\phi} (I^s_{eji\phi})^H$. For $i \in \mc N$ and $\phi \in \mc P$, $\bs S^g_{k\phi} = V_{i\phi} (I_{k\phi}^g)^H, \ \forall k \in \mc G(i)$ and $\bs S^b_{l\phi} = V_{i\phi}(I^b_{l\phi})^H, \ \forall l \in \mc L(i)$. 
\end{itemize}

We introduce auxiliary variables $\bs W_i$ for $\bs V_i \bs V_i^H$, $\forall i \in \mathcal{N}$, $\bs{M}_{eij}$ for $\bs{V}_i (\bs{I}_{eij}^s)^H$ and $\bs{L}_{eij}$ for $\bs{I}_{eij}^s(\bs{I}_{eij}^s)^H$, $\forall (e,i,j) \in \mc E$, using \eqref{eq:acpf:rank1:leafbus} and \eqref{eq:acpf:rank1}. We can replace Kirchhoff's current laws for buses and lines in (ii) with \eqref{eq:acpf:power-bal} and \eqref{eq:acpf:power-loss} by taking the conjugate transpose of both sides, multiplying both sides by $\mathbf{V}_i$ from the left, and applying Ohm's law and the power equations given in (i) and (iii) respectively. In summary, these physical laws can be represented as the following set of constraints, often referred to as a branch-flow model \cite{gan2014convex}:
\begin{subequations}
\begin{align}
  & \bs W_i = \bs V_i \bs V_i^H, \forall i \in \mathcal N_{{leaf}},  \label{eq:acpf:rank1:leafbus}\\
  & \forall (e,i,j) \in \mc E, \nonumber\\
  & \quad \left[\begin{array}{cc}
    \bs W_i & \bs M_{eij}\\
    \bs M^H_{eij} & \bs L_{eij}
    \end{array}\right] = \left[\begin{array}{c}
      \bs V_i \\
      \bs I^s_{eij}
      \end{array}\right] \left[\begin{array}{c}
        \bs V_i \\
        \bs I^s_{eij}
        \end{array}\right]^H,\label{eq:acpf:rank1}\\
  & \quad \bs V_j = \bs V_i - {\color{black}\bs Z_{e}} \bs I^s_{eij}, \label{eq:acpf:ohm}\\
  & \quad I^s_{eij\phi} = 0, \ \forall \phi \notin \mc{P}_e, \label{eq:acpf:current_on_missing_phase}\\
  & \quad \bs S_{eij} = \mbox{diag}(\bs M_{eij}), \label{eq:acpf:power-eq}\\
  & \quad \bs S_{eji} = - \bs S_{eij} + \mbox{diag}( {\color{black}\bs Z_{e}} \bs L_{eij}),\label{eq:acpf:power-loss}\\
  & \bs S^{b}_{l} = {\color{black} \bs S^0_l}, \forall l \in \mc L,  \label{eq:acpf:load-power}\\
  & \forall i \in \mc N, \nonumber\\
  & \quad \sum_{(e,i,j) \mbox{ or } (e,j,i) \in \mc E} (\bs S_{eij} + \mbox{diag}(\bs W_i({\color{black}\bs Y^{sh}_{eij}})^H)) =\nonumber \\
  & \quad -\sum_{s \in \mathcal B(i)} \mbox{diag}(\bs W_i  ({\color{black}\bs Y^{sh}_{s}})^H) + \sum_{k \in \mathcal G(i)} \bs S^{g}_{k} - \sum_{l \in \mathcal L(i)} \bs S^{b}_{l}.\label{eq:acpf:power-bal}
\end{align}
\label{eqs:acpf}
\end{subequations}

\subsection{Multiphase AC OPF}
Each $k \in \mc G$ incurs cost $C_k(\bs S^g_k)$ for generating $\bs S^g_k \in \mathcal S_k$. OPF aims to minimize the total generation cost for satisfying the load while abiding by the physical laws, namely \eqref{eqs:acpf}, as well as a set of operational limits. These constraints include generation bounds, 
voltage bounds, and substation voltage: 
\begin{subequations}
  \begin{align}
    & \bs S^g_k \in \mc S_k, \ \forall k \in \mc G, \label{eq:op-constr:gen-bounds}\\
    & {\color{black}\underline v_{i \phi}^2} \le \mbox{diag}(\bs W_{i})_\phi \le {\color{black}\overline v_{i\phi}^2}, \ \forall i \in \mc N, \ \phi \in \mc P,  \label{eq:op-constr:volt-bounds}\\
    & \bs W_0 = {\color{black}\bs V_0^{ref}(\bs V_0^{ref})^H}. \label{eq:op-constr:ref-volt}
  \end{align}
  \label{eqs:op-constr}
  \end{subequations}
  In summary, the OPF can be formulated as
  \begin{equation*}
    \mbox{(AC)}: \min \left\{\sum_{k \in \mc G} C_k(\bs S^g_k): \eqref{eqs:acpf}, \eqref{eqs:op-constr}\right\}.
    \label{prob:acopf}
  \end{equation*}

\subsection{Linear approximation of (AC)}\label{sec:linearized}
Gan et al.~\cite{gan2014convex} proposed a linear approximation of (AC) based on the assumption below:
\begin{assumption}
  \begin{enumerate}
    \item[(a)] Line losses are small; that is,  $\mbox{diag}({\color{black}\bs Z_e} \bs L_{eij}) \ll \bs S_{eij}$ for $(e,i,j) \in \mc E$.
    \item[(b)] Voltages are nearly balanced; for example,  if $\mathcal P_i = \{a,b,c\}$, then 
    $$\frac{V_{ia}}{V_{ib}} \approx \frac{V_{ib}}{V_{ic}} \approx \frac{V_{ic}}{V_{ia}} \approx e^{j 2 \pi / 3}.$$
  \end{enumerate}\label{assum:linearized}
\end{assumption}

By Assumption \ref{assum:linearized} (b), for each $(e,i,j) \in \mc E$, $\bs M_{eij}$ can be approximated with $\bs S_{eij}$ as follows:
\begin{equation}\bs M_{eij} = \bs V_i (\bs I_{eij}^s)^H \approx \bs\Gamma \mbox{diag} (\bs S_{eij}), \label{eq:lp:M}
\end{equation}
where $\gamma = e^{-\mathrm{i} 2\pi/3}$ and $\bs\Gamma = \left[\begin{array}{ccc}
  1 & \gamma^2 & \gamma \\ 
  \gamma & 1 & \gamma^2\\
  \gamma^2 & \gamma & 1
\end{array}\right]$. 

In addition, by multiplying both sides of \eqref{eq:acpf:ohm} and \eqref{eq:acpf:current_on_missing_phase} by their conjugate transposes, we can further eliminate $\bs I^s_{eij}$ and $\bs V_i$ and obtain
\begin{subequations}
\begin{align}
  &\bs W_j = \bs W_i + {\color{black}\bs Z_e} \bs L_{eij}{\color{black} \bs Z_e^H} - \bs M_{eij}  {\color{black}\bs Z_e^H}  - {\color{black}\bs Z_e} \bs M_{eij}^H,\label{eq:sdp:ohm}\\
  &\bs L_{eij\phi\cdot} = 0, \ \forall \phi \notin \mc P_e,\label{eq:sdp:current_on_missing_phase}
\end{align}
\label{eqs:sdp:ohm}
\end{subequations}
\noindent where $\bs L_{eij\phi\cdot}$ denotes the row of $\bs L_{eij}$ corresponding to $\phi$. By Assumption \ref{assum:linearized} (a), terms with ${\color{black}\bs Z_e} \bs L_{eij}$ in \eqref{eq:sdp:ohm} and \eqref{eq:acpf:power-loss} can be neglected, and we have
\begin{subequations}
  \begin{align}
    &\bs W_j = \bs W_i - \bs M_{eij}  {\color{black}\bs Z_e^H}  - {\color{black}\bs Z_e} \bs M_{eij}^H,\label{eq:lp:ohm}\\
    &\bs S_{eji} = - \bs S_{eij},\label{eq:lp:power-loss}
  \end{align}\label{eq:lp:ex-loss}
\end{subequations}
\noindent and \eqref{eq:sdp:current_on_missing_phase} can be replaced by 
\begin{equation}
  \bs S_{eij\phi} = 0, \ \forall \phi \notin \mc P_e.
  \label{eq:lp:current_on_missing_phase}
\end{equation} Thus, we can exclude $\bs L_{eij}$ from the model. 
As a result, we obtain the following linear approximation:
\begin{subequations}
  \begin{align}
    \mbox{(LP)}:  \min \ & \sum_{k \in \mc G} C_k(\bs S^g_k) \nonumber\\
    \mbox{s.t.} \ & \eqref{eq:lp:M}, \eqref{eq:lp:ex-loss}, \eqref{eq:lp:current_on_missing_phase}, \forall (e,i,j) \in \mc E, \label{eq:acopf:lp:1}\\
    &\eqref{eq:acpf:load-power}-\eqref{eq:acpf:power-bal}, \eqref{eqs:op-constr}.\label{eq:acopf:lp:2}
  \end{align}  \label{eqs:acopf:lp}
\end{subequations}

\section{Exponential load with delta connection} \label{sec:vol-dependent-load-model}
Most multiphase OPF models, including the OPF models illustrated in Section \ref{sec:prelim}, assume that the amount of power withdrawn from bus $i$ by load $l \in \mc L(i)$ (i.e., $\bs S^b_l$) is constant. In practice, however, $\bs S^b_l$ may change depending on the applied voltage as well as how the load is connected to the bus (e.g., delta or wye). 
In order to better incorporate the load behavior, a more general load model for voltage-dependent loads configured either as wye or delta can be included in the OPF models. 

\subsection{Configuration of multiphase load}\label{sec:linearized-model:config}
\begin{figure}[!t]  
  \centering
  \begin{subfigure}[b]{0.23\textwidth}
      \centering      
      \includegraphics[width=\textwidth]{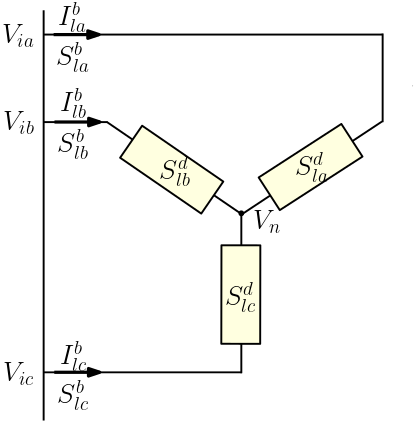}                       
      \caption{Wye connection}\label{fig:load:connection:wye}
  \end{subfigure}
  \begin{subfigure}[b]{0.21\textwidth}
      \centering      
      \includegraphics[width=\textwidth]{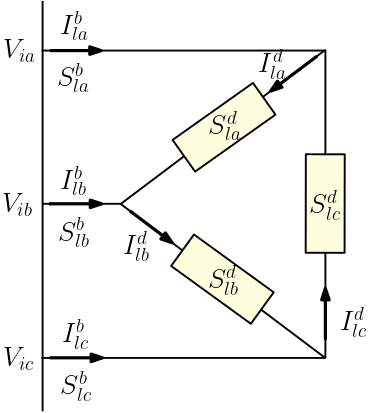}                       
      \caption{Delta connection}\label{fig:load:connection:delta}
  \end{subfigure} 
  \caption{Three-phase load configuration} 
  \label{fig:load:connection}
\end{figure}

As illustrated in Figure \ref{fig:load:connection}, depending on how each individual load of a multiphase load is connected, the relationship between $\bs S^b_l$ and $\bs S^d_l=[S^d_{l\phi}]_{\phi \in \mc P}$ as well as the voltage applied to each individual load differs. Let $\mc Y \subseteq \mc L$ and $\mc D \subseteq \mc L$ respectively denote the set of wye- and delta-connected loads. Throughout this paper we assume that the neutral conductor is perfectly grounded (i.e., $V_n=0$).

Figure \ref{fig:load:connection}(a) illustrates a three-phase wye-connected load $l \in \mc Y$ at bus $i \in \mc N$. Note that for each $\phi \in \mathcal P$, the voltage applied to the load (i.e., $V_{l\phi}$) is $V_{i\phi}$ and the current passing through the load (i.e., $I^d_{l\phi}$) is $I_{l\phi}^b$ and thus 
\begin{equation}
  \bs S^{b}_{l} = \bs S^d_{l}.
  \label{eq:wye:power}
\end{equation}

For a three-phase delta-connected load $l \in \mc D$ at bus $i \in \mc N$, as illustrated in Figure \ref{fig:load:connection}(b), the voltage applied to the load is 
\begin{equation}
\label{eq:delta_V}
\bs V_{l} = {\bs \Lambda} \bs V_{i}, 
\end{equation}
where 
$${\bs \Lambda} = \left[\begin{array}{ccc}
  1 & -1 & 0\\
  0 & 1 & -1\\
  -1 & 0 & 1
\end{array}\right].$$

In addition, from Kirchhoff's current law, we have 
\begin{equation}
  \bs I^b_{l} = {\bs \Lambda}^{\top} \bs I_l^d. 
  \label{eq:Kirchhoff}
\end{equation}
Therefore, 
  \begin{subequations}
    \begin{align}
  &\bs S^d_{l} = \mbox{diag}(\bs V_l (\bs I_{l}^d)^H) = \mbox{diag}({\bs \Lambda} \bs V_{i} (\bs I_{l}^d)^H),\\
  &\bs S^{b}_{l}= \mbox{diag}(\bs V_i (\bs I^b_l)^H)= \mbox{diag}(\bs V_i (\bs I^d_l)^H{\bs \Lambda}).
    \end{align}\label{eqs:delta:power}
\end{subequations}
By replacing each $\bs I^b_l$ with $\bs I^d_l$ for $l \in \mc Y$ and with ${\bs \Lambda}^\top \bs I^d_l$ for $l \in \mc D$ in the physical laws and introducing auxiliary variables $\bs L^d_{eij}$ for $\bs I^d_{l} (\bs I^d_{l})^H$ and $\bs X_{l}$ for $\bs V_i (\bs I^d_{l})^H$, we can derive the AC OPF problem with delta connections as follows:
\begin{subequations}
\begin{align}
  \mbox{(AC-D): }\min \ & \sum_{k \in \mc G} C_k(\bs S^g_k)\nonumber\\
  \mbox{s.t.}\ & \eqref{eq:acpf:rank1}-\eqref{eq:acpf:power-loss},\eqref{eq:acpf:power-bal}, \eqref{eqs:op-constr},\label{eq:acpf-d:acpf}\\
  & \bs S_l^d = {\color{black} \bs S^0_l}, \forall l \in \mc L, \label{eq:acpf-d:load-power}\\
  & \forall i \in \mc N,  l \in \mc D(i), \nonumber\\
  & \quad \left[\begin{array}{cc}
    \bs W_i & \bs X_{l}\\
    \bs X_{l}^H & \bs L^d_{l}
    \end{array}\right] = \left[\begin{array}{c}
      \bs V_i \\
      \bs I^d_{l}
      \end{array}\right] \left[\begin{array}{c}
        \bs V_i \\
        \bs I^d_{l}
        \end{array}\right]^H,\label{eq:acpf-d:rank1}\\
  & \quad  \bs S_l^d = \mbox{diag}(\bs \Lambda \bs X_l), \label{eq:acpf-d:delta:power-connection1}\\
  & \quad  \bs S_l^b = \mbox{diag}(\bs X_l \bs \Lambda), \label{eq:acpf-d:delta:power-connection2}\\
  & \quad  I^d_{l\phi} = 0, \ \forall \phi \notin \mc P_l,\label{eq:acpf-d:delta:current-on-missing-phase}\\
  & \eqref{eq:wye:power}, \forall i \in \mc N,  l \in \mc Y(i). \label{eq:acpf-d:wye:power-connection}
\end{align}\label{eqs:acopf-d}
\end{subequations}
\subsection{Linear approximation of (AC-D)}
We propose a linear approximation of \eqref{eqs:delta:power} that can be readily added to (LP) under Assumption \ref{assum:linearized}. 
Note that by taking the conjugate transpose of each side of \eqref{eq:Kirchhoff} and multiplying  both sides with $\bs V_i$ from the left, we have
$\sum_{\phi \in \mc P}S^b_{l\phi} = \sum_{\phi \in \mc P} {S^d_{l\phi}}$. In addition, by Assumption \ref{assum:linearized} (b), it follows that 
\begin{equation}
  {S_{la}^d} = (V_{ia} - \gamma V_{ia})(I^d_{la})^H  \Longrightarrow V_{ia} (I^d_{la})^H = \frac{1}{1-\gamma} {S_{la}^d},
  \label{eq:delta:Sd1}
\end{equation}
where $\gamma = e^{-\mathrm{i} (2\pi/3)} = - \frac{1}{2} - \mathrm{i} \frac{\sqrt{3}}{2}$. Using \eqref{eq:delta:Sd1} and Assumption \ref{assum:linearized} (b), we can express $S^d_{lb}$ and $S^d_{lc}$ as follows:
  \begin{align*}
    {S_{lb}^d} &= (\gamma V_{ia} - \gamma^2 V_{ia}) (I_{lb}^d)^H  = (1-\gamma)\gamma V_{ia} (I_{lb}^b + I_{la}^d)^H\\
    & = (1-\gamma) \gamma \left[\frac{1}{\gamma} S^b_{lb} + \frac{1}{1-\gamma} {S^d_{la}}\right],\\
    {S_{lc}^d} &= (\gamma^2 V_{ia}-V_{ia}) (I_{lc}^d)^H = (\gamma^2 -1) V_{ia}(I_{lb}^b + I_{lc}^b + I_{la}^d)^H\\
    & = (\gamma^2 - 1) \left[\frac{1}{\gamma} S^b_{lb} + \frac{1}{\gamma^2} S^b_{lc} + \frac{1}{1-\gamma} {S^d_{la}}\right].
  \end{align*}  

By expressing the equations in real and imaginary parts, we derive the following linear system of equations that connect $\bs S^d_l$ and $\bs S^b_l$ for each delta load $l \in \mathcal D$:
\begin{subequations}
  \begin{align}
    &p_{l1}^b + p_{l2}^b + p_{l3}^b = {p_{l1}^d + p_{l2}^d + p_{l3}^d},\\
    &q_{l1}^b + q_{l2}^b + q_{l3}^b = {q_{l1}^d + q_{l2}^d + q_{l3}^d},\\
    & \frac{3}{2} p^{b}_{l2} - \frac{\sqrt{3}}{2} q^{b}_{l2}  = {p^{d}_{l2} + \frac{1}{2} p^{d}_{l1} - \frac{\sqrt{3}}{2} q^{d}_{l1}},\\
    & \frac{\sqrt{3}}{2} p^{b}_{l2} + \frac{3}{2} q^{b}_{l2} = {\frac{\sqrt{3}}{2} p^{d}_{l1} +\frac{1}{2} q^{d}_{l1} + q^{d}_{l2}},\\
    & \sqrt{3} q^{b}_{l2} + \frac{3}{2} p^{b}_{l3}-\frac{\sqrt{3}}{2} q^{b}_{l3} = {\frac{1}{2} p^{d}_{l1} + \frac{\sqrt{3}}{2} q^{d}_{l1} + p^{d}_{l3}}, \\
    & - \sqrt{3}p^{b}_{l2} + \frac{\sqrt{3}}{2} p^{b}_{l3} + \frac{3}{2} q^{b}_{l3}   = {- \frac{\sqrt{3}}{2} p^{d}_{l1} + \frac{1}{2} q^{d}_{l1}  + q^{d}_{l3}}.
  \end{align}  \label{eq:delta:power}
\end{subequations}
\begin{remark}
  Note that \eqref{eq:delta:power} is a system of 6 linear equations with 6 unknowns, say, $A\bs x = b$. Since the corresponding square matrix $A \in \mathbb R^{6 \times 6}$ is invertible, $\bs p^b_l$ and $\bs q^b_l$ are uniquely defined once $\bs p^d_l$ and $\bs q^d_l$ are determined. This implies that under Assumption \ref{assum:linearized}, \eqref{eq:delta:power} exactly describes the relationship between $\bs S^b$ and $\bs S^d$ for delta-connected devices.
\end{remark}

To summarize, the linear approximation of (AC-D) can be expressed as follows:
\begin{subequations}
\begin{align}
  \mbox{(LP-D)}:  \min \ & \sum_{k \in \mc G} C_k(\bs S^g_k)\nonumber\\
  \mbox{s.t.} \ &\eqref{eqs:acopf:lp}, \eqref{eq:acpf-d:load-power}, \\
  & \eqref{eq:delta:power}, \forall l \in \mc D,\label{eq:acopf-d:lp:delta}\\
  & \eqref{eq:wye:power}, \forall l \in \mc Y. \label{eq:acopf-d:lp:wye}
\end{align}
\label{eqs:acopf-d:lp}
\end{subequations}
\begin{remark}
  Note that the linear approximation of the delta connection (i.e., \eqref{eq:delta:power}) can also be applied to any delta-connected network components (e.g., generators in delta connections). It approximates the bus injection/withdrawal $\bs S^b$ of a delta-connected device based on its power production/consumption $\bs S^d$.
\end{remark}

\subsection{Exponential load model}\label{sec:exponential-load}
In the preceding sections we assumed that the power consumption of load $l \in \mc L$ at phase $\phi$ (i.e., $S^d_{l\phi}$) is constant at its nominal power $\color{black}S^0_{l\phi} = p^0_{l\phi} +\mathrm{i}q^0_{l\phi}$. In practice, however, $\bs S^d_l$ changes depending on the applied voltage. A widely used model for voltage-dependent load is an exponential load model that assumes that the power consumption of a load is proportional to the applied voltage magnitude raised to some power \cite{korunovic2018recommended}. 
For each multiphase load $l \in \mc L$, its power consumption at phase $\phi \in \mathcal P$, namely., $S^d_{l\phi} = p^d_{l\phi} + \mathbf{i} q^d_{l\phi}$, is computed by
\begin{align}
  p^d_{l\phi} = {\color{black}p^{0}_{l\phi}} \left(\frac{|V_{l\phi}|}{{|\color{black}V^0_{l\phi}|}}\right)^{{\color{black}\alpha_{l\phi}}}, \ \ q^d_{l\phi} = {\color{black}q^{0}_{l\phi}} \left(\frac{|V_{l\phi}|}{{|\color{black}V^0_{l\phi}|}}\right)^{{\color{black}\beta_{l\phi}}}, \label{volt-dependent-model}
\end{align}
where 
$p^d_{l\phi}$ (resp., $q^d_{l\phi}$) denotes the active (resp., reactive) power consumption of $l$ at phase $\phi$ when the magnitude of voltage applied to the load is $|V_{l\phi}|$ and where
${\color{black}p^{0}_{l\phi}}$, ${\color{black}q^{0}_{l\phi}}$, and ${|\color{black}V^0_{l\phi}|}$ respectively denote reference values for active power load, reactive power load, and voltage magnitude applied to the load, which are given as data. The exponents ${\color{black}\alpha_{l\phi}}$ and ${\color{black}\beta_{l\phi}}$ are also given as data, which are nonnegative numbers that characterize the voltage dependency of load $l$. For instance, some special choice of the exponents yields classical load models; With the exponents ${\color{black}\alpha_{l\phi}}$ and ${\color{black}\beta_{l\phi}}$ equal to $0$, $1$, and $2$, \eqref{volt-dependent-model} represents constant power, constant current, and constant impedance load, respectively. Exponent values other than 0, 1, or 2 can be employed to model more general load types; see, for example, \cite{korunovic2018recommended}.

We denote the coefficients of \eqref{volt-dependent-model} by $\color{black}a_{l\phi}$ and $\color{black}b_{l \phi}$ (i.e., $\color{black}a_{l\phi} = \frac{p_{l\phi}^{0}}{|V^0_{l\phi}|^{\alpha_{l\phi}}}$ and $\color{black}b_{l\phi} = \frac{q_{l\phi}^{0}}{|V^0_{l\phi}|^{\beta_{l\phi}}}$). Without loss of generality, we assume that $\color{black}a_{l\phi}$ and $\color{black}b_{l\phi}$ are nonzeros, since otherwise we can fix the corresponding values of $p^d_{l\phi}$ and $q^d_{l\phi}$ at zeros. Now we may express \eqref{volt-dependent-model} in terms of the squared magnitude of voltage applied to load $l$ at phase $\phi$, denoted by $v_{l\phi} = |V_{l\phi}|^2$:
\begin{align}
  \frac{1}{{\color{black}a_{l\phi}}}p^d_{l\phi} =  v_{l\phi}^{\frac{\color{black}\alpha_{l\phi}}{2}}, \ \ \frac{1}{{\color{black}b_{l\phi}}}q^d_{l\phi} =  v_{l\phi}^{\frac{\color{black}\beta_{l\phi}}{2}},\label{volt-dependent-model-2}
\end{align}
where 
\begin{subequations}
  \begin{align}
    &\bs v_{l} =\mbox{diag}(\bs V_i \bs V_i^H) = \mbox{diag}( \bs W_i), \ \forall l \in \mc Y(i),\label{eq:exp-load:volt:wye}  \\
    &\bs v_{l} = \mbox{diag}(\bs \Lambda\bs V_i (\bs \Lambda\bs V_i)^H) =\mbox{diag}(\bs \Lambda \bs W_i \bs \Lambda^\top), \ \forall l \in \mc D(i).
  \end{align}\label{eqs:exp-load:volt}  
\end{subequations}
To summarize, AC OPF with delta-connected exponential loads can be posed as follows:
  \begin{align*}
    \mbox{(AC-D-E)}: \min \ & \sum_{k \in \mc G} C_k(\bs S^g_k)\\
    \mbox{s.t.} \ & \eqref{eq:acpf-d:acpf}, \eqref{eq:acpf-d:rank1}-\eqref{eq:acpf-d:wye:power-connection}, \eqref{volt-dependent-model-2}, \eqref{eqs:exp-load:volt}.
  \end{align*}

\subsection{Linear approximation of (AC-D-E)}
We utilize the linear approximation of \eqref{volt-dependent-model-2} at $v_{l\phi} = 1$, implemented in \cite{fobes2020powermodelsdistribution}:
\begin{subequations}
\begin{align}
  p^d_{l\phi} = \frac{\color{black}a_{l\phi} \alpha_{l\phi}}{2} (v_{l\phi}-1)+{\color{black}a_{l\phi}}, \\ 
  q^d_{l\phi} = \frac{\color{black}b_{l\phi} \beta_{l\phi}}{2} (v_{l\phi}-1) + {\color{black}b_{l\phi}}. 
\end{align}\label{eq:volt-dependent-model-linearized}
\end{subequations}
Note that when the exponents are 0 or 2, the linearized equation correctly represents the constant power and constant impedance load models, respectively.

 For delta-connected loads, the voltage applied to the load is $V_{i\phi}-V_{i\phi^+}$ for each $\phi \in \mathcal P$, and thus $v_{l\phi} = |V_{i\phi}-V_{i\phi^+}|^2$. Using the assumption that voltages are almost balanced (i.e., Assumption \ref{assum:linearized} (b)), we approximate $V_{i\phi+}$ with $\gamma V_{i\phi}$, and thus 
\begin{equation}
  v_{l\phi} = |1-\gamma|^2 W_{i\phi\phi} = 3W_{i\phi\phi}.
  \label{eq:linear:delta:volt}
\end{equation}

In summary, the linear approximation of (AC-D-E) can be expressed as follows:
  \begin{align*}
    \mbox{(LP-D-E)}: \min \ & \sum_{k \in \mc G} C_k(\bs S^g_k)\\
    \mbox{s.t.} \ & \eqref{eqs:acopf:lp},\eqref{eq:acopf-d:lp:delta}, \eqref{eq:acopf-d:lp:wye},\eqref{eq:volt-dependent-model-linearized}, \eqref{eq:exp-load:volt:wye}, \eqref{eq:linear:delta:volt}.
  \end{align*}

\section{Numerical results}\label{sec:experiments}
In this section we analyze the performance of the proposed linear approximation of (AC-D-E), namely (LP-D-E), on the IEEE 13, 37, and 123 bus systems. The generation cost function $C_k(\bs S_k^g)$ for each $k \in \mc G$ is defined to be the sum of its real power generation $\sum_{\phi \in \mc P_k} p_{k\phi}$. In these test systems, only the voltage source serves as a generator, and thus the problem minimizes the total import from the transmission grid. The voltage bounds are set to be [0.8 p.u., 1.2 p.u.]. Substation transformers and regulators are removed, and the load transformers are modeled as lines with equivalent impedance. 
We eliminate switches or regard them as lines according to their default status. All capacitors are assumed active and generate reactive power according to their ratings and the associated bus voltages. Substations voltages are set as $V_0^{ref} = \bar V[1, e^{-j2\pi/3},e^{j2\pi/3}]^T$, where $\bar V$ is given as the substation voltage in the test systems. This experiment setup is commonly used in the literature \cite{zhao2017optimal,zhou2021exactness}. 

\begin{figure}
\centering
  \includegraphics[width=0.9\linewidth]{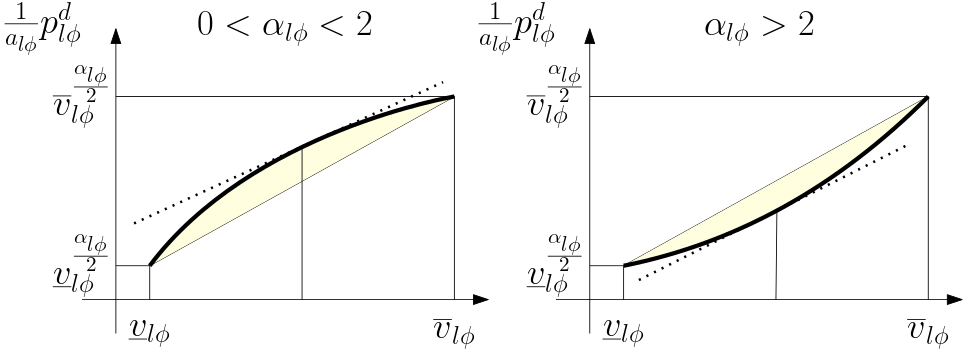}
  \caption{Graphical illustration of \eqref{volt-dependent-model-2} for $\alpha_{l\phi} \neq 0, 2$ (thickened curves), its convex relaxation proposed in \cite{claeys2021voltage} (shaded regions), and its linear approximation given in \eqref{eq:volt-dependent-model-linearized} (dotted lines)}
  \label{fig:load-convex}
\end{figure}
In addition to comparing the results of (AC-D) and (LP-D), as well as the results of (AC-D-E) and (LP-D-E), we also make comparisons with the following two convex relaxations:
\begin{itemize}
    \item (SDP-D): a SDP relaxation of (AC-D), proposed in \cite{zhou2021exactness}. 
    \item (CONE-D-E): a convex relaxation of (AC-D-E) proposed in \cite{claeys2021voltage}. It adds to (SDP-D) a convex conic relaxation of \eqref{volt-dependent-model-2}, which is illustrated in Figure \ref{fig:load-convex}. 
\end{itemize} 
For both (SDP-D) and (CONE-D-E), we set the penalty weight $\rho$ for minimizing the rank of $\bs L_l^d$ to be 100. Once they are solved, the resultant penalty term is subtracted from the objective value to calculate the true objective value. 

All experiments were executed on a machine with 32 GB of memory on an Intel Core i7 at 2.3 GHz. (AC-D) and (AC-D-E) were solved by using Ipopt 3.14.4, running with linear solver MUMPS 5.4.1, (SDP-D) and (CONE-D-E) were by Mosek 9.3, and the linear models were via CPLEX 20.1.

\subsection{Performance of (SDP-D) and 
(LP-D) solutions}\label{sec:exp:sdp-d-exactness}

To see the performance of the linearized model for the delta connection, we first compare (AC-D), (SDP-D), and (LP-D) by fixing all loads at their nominal power (i.e., the constant power model). As noted in \cite{zhou2021exactness}, we observe that (SDP-D) produces almost rank-1 solutions for all the test systems, 
implying the global optimality of the (SDP-D) solutions. We also observe that (AC-D) produces globally optimal solutions that comply with (SDP-D), providing the same objective values, $\bs w$, and $\bs S^b$ as (SDP-D) for all the test instances.

\begin{table}[t!]
  \caption{Performance of (LP-D) with respect to (SDP-D)}
  \label{table:lp-d}
  \centering
  \begin{tabular}{ccccccccc}
    \toprule
    \multirow{2}{*}{Instance} &\multicolumn{3}{c}{LP-D}\\
    \cmidrule(lr){2-4}
     & $\Delta\bs w(\%)$ &  $\Delta \bs p^b (\%) $ &
     $\Delta \bs q^b (\%) $ \\  
    \midrule
    IEEE 13  & 0.93 & 0.55 & 3.32\\
    IEEE 37  & 0.12 & 0.5 & 2.1\\
    IEEE 123 & 0.41 & 0.07 & 0.59\\
    \bottomrule
  \end{tabular}
\end{table}
\begin{table}[t!]
  \caption{Solution time in seconds}
  \label{table:computational}
  \centering
  \begin{tabular}{ccccccccc}
    \toprule
    Instance & AC-D & SDP-D & LP-D\\
    \midrule
    IEEE 13  & 2.31 & 0.92 & 0.003\\
    IEEE 37  & 3.51 & 2.50 & 0.006\\
    IEEE 123 & 4.74 & 6.55 & 0.059\\
    \bottomrule
  \end{tabular}
\end{table}

Table \ref{table:lp-d} compares the solutions of (LP-D) with those of (SDP-D) to see the performance of the linearization, where $\Delta \bs x$ represents the average relative difference in variable $\bs x$ in percentage; that is,  $\Delta \bs x = 100 \times \sum_{i=1}^{n}\frac{|x_i^{LP} - x_i^{SDP}|}{|x_i^{SDP}|}/n$, where $n$ denotes the dimension of $\bs x$. The result suggests that the error does not exceed 1\% for the voltage magnitudes and the real power withdrawals for all the instances. 

\subsection{Performance of (CONE-D-E) and (LP-D-E)  solutions}
\begin{table}[t!]
  \caption{Objective value comparison}
  \label{table:objval-e}
  \centering
  {\color{black}
  \begin{tabular}{ccccccccc}
    \toprule
    Instance & AC-D-E & CONE-D-E & LP-D-E\\
    \midrule
    IEEE 13  & 3520.35 & 3461.98 & 3421.81 \\
    IEEE 37  & 2478.05 & 2325.19 & 2480.39 \\
    IEEE 123 & 3496.30 & 3412.49 & 3430.68 \\
    \bottomrule
  \end{tabular}
  }
\end{table}
\begin{table}[t!]
  \caption{Relative difference in percentage from solutions of (AC-D-E)}
  \label{table:sol-e}
  \centering
  {\color{black}
  \begin{tabular}{cccccccccc}
    \toprule
    \multirow{2}{*}{Instance} & \multicolumn{3}{c}{CONE-D-E} & \multicolumn{3}{c}{LP-D-E}\\
    \cmidrule(lr){2-4} \cmidrule(lr){5-7}
     & $\Delta \bs w$ & $\Delta \bs p^b$ & $\Delta \bs q^b$ & $\Delta \bs w$  & $\Delta \bs p^b$ & $\Delta \bs q^b$ \\  
    \midrule
    IEEE 13  & 0.15 & 2.01 & 1.82 & 0.74 & 0.78 & 3.8\\
    IEEE 37  & 0.33 & 5.81 & 15.99 & 0.06 & 2.98 & 8.35\\
    IEEE 123 & 0.15 & 0.93 & 1.46 & 0.3 & 0.38 & 0.7\\
    \bottomrule
  \end{tabular}
  }
\end{table}
\begin{table}[t!]
  \caption{Solution time in seconds}
  \label{table:computational-e}
  \centering
  \begin{tabular}{ccccccccc}
    \toprule
    Instance & AC-D-E & CONE-D-E & LP-D-E\\
    \midrule
    IEEE 13  & 2.38 & 1.60 & 0.004\\
    IEEE 37  & 2.93 & 2.98 & 0.007\\
    IEEE 123 & 5.61 & 5.98 & 0.078\\
    \bottomrule
  \end{tabular}
\end{table}
In this section we analyze the performance of the formulations with exponential load models. 
We observe that (CONE-D-E) returns rank-1 solutions for all the instances; however, it no longer provides globally optimal solutions to (AC-D-E) and can only serve as a lower bound since the conic reformulation of the exponential load model may not be exact unless all loads have $\alpha = 0 \mbox{ or }2$ and $\beta = 0 \mbox{ or } 2$. Since all the test cases include loads with exponents that are not 0 or 2, (CONE-D-E) provides lower bounds, as shown in Table \ref{table:objval-e}. Table \ref{table:sol-e} summarizes the average relative differences between solutions of (AC-D-E) and those of (CONE-D-E) and (LP-D-E). The relative difference between the LP solution and the AC solution does not exceed 1\% with regard to the voltage magnitudes and 3\% for the real power withdrawals. 

\begin{figure}[t!]
  \centering
  \includegraphics[width=0.6\linewidth]{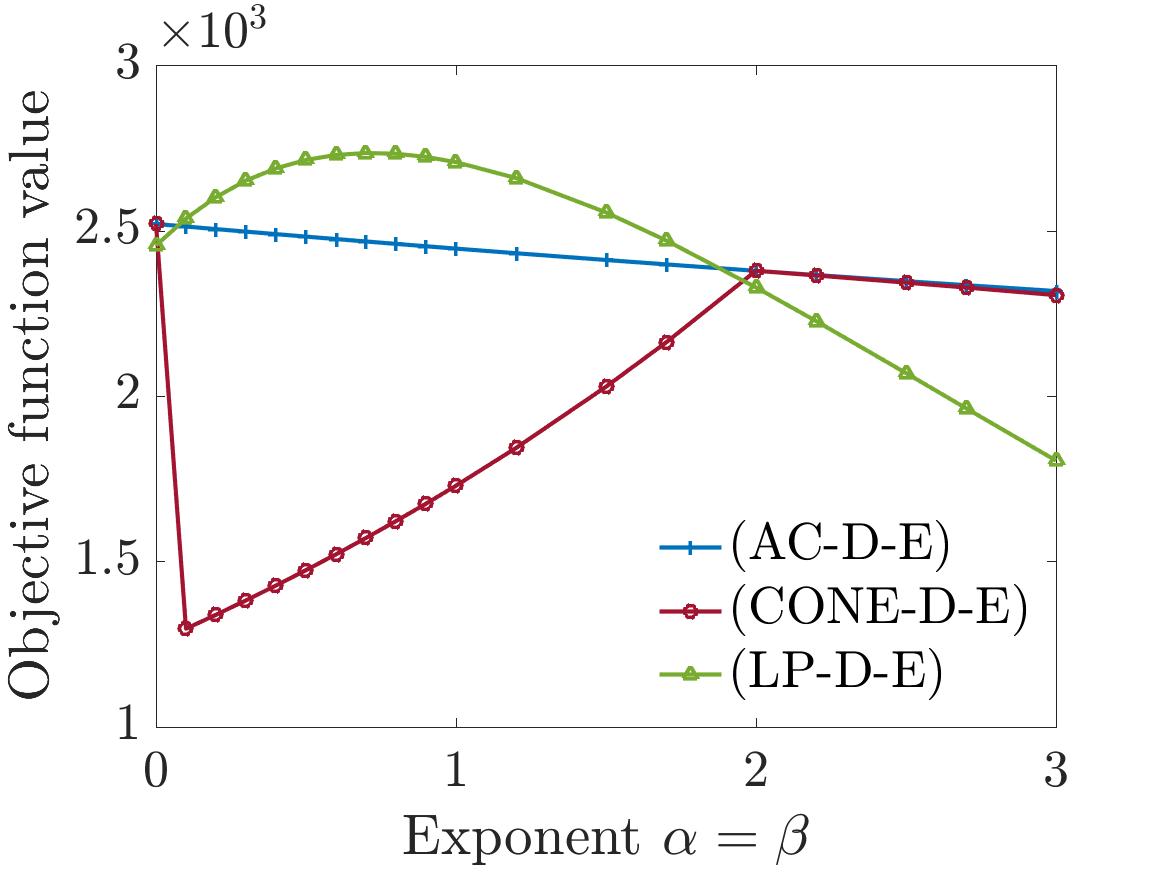}
  \caption{{\color{black}Comparison of (AC-D-E), (CONE-D-E), and (LP-D-E)}}
  \label{fig:delta-exp}
\end{figure}
{\color{black}Furthermore, to evaluate the performance of (CONE-D-E) and (LP-D-E) under varying exponents of the exponential model, we conduct experiments on the IEEE 37-bus system, treating all loads as exponential. The experiment involves uniformly varying the exponents of all loads (i.e., $\alpha_{l,\phi}, \beta_{l,\phi} \ \forall l \in \mc L, \phi \in \mc P_l$) from 0 to 3. Figure \ref{fig:delta-exp} plots the objective function values of (AC-D-E), (CONE-D-E), and (LP-D-E) for different choices of the exponent. For exponents $0 \leq \alpha=\beta < 2$, (LP-D-E) tends to produce objective function values closer to (AC-D-E), while for $\alpha=\beta \geq 2$, (CONE-D-E) outperforms (LP-D-E). This behavior is expected from the characteristics of the convex relaxation and the linear approximation of the exponential load model, depicted in Figure \ref{fig:load-convex}. Note that for $0 \le \alpha = \beta \le 2$, the linear approximation overestimates the total power consumption, and otherwise, it is an underestimation. This explains why (LP-D-E) has a lower objective value for $0 \le \alpha = \beta \le 2$ and a higher objective value for $\alpha = \beta > 2$ than those of (AC-D-E). When $\alpha=\beta=$ 0 or 2 (i.e., when the linearized load model is exact), the objective function value of (LP-D-E) is slightly smaller than that of (AC-D-E) since it ignores power losses. On the contrary, for the convex relaxation, since the objective is to minimize the total power generation, $p^d_{l\phi}$ is more likely to occur on the linear underestimator 
when $0 \le \alpha = \beta \le 2$ and on the curve $v_{l\phi}^{\alpha_{l\phi}/2}$ for $\alpha = \beta > 2$ (see, e.g., Figure \ref{fig:load-convex}). 
This explains why (CONE-D-E) performs poorly for $0 \le \alpha = \beta \le 2$ and becomes accurate for $\alpha = \beta > 2$. Hence, depending on the exponent of each voltage-dependent load, one might choose either the linear approximation or the conic relaxation on a load-by-load basis to enhance accuracy.}

\subsection{Solution Time}
Tables \ref{table:computational} and \ref{table:computational-e} provide a comparison of the time required to solve each model. The proposed linear approximation yields solutions faster than the benchmarked conic relaxations by at least two orders of magnitude.

\subsection{Importance of modeling delta-connected, voltage-dependent loads}
\begin{figure*}[t]
  \centering
  \subfloat[][(AC-W-E)]{\includegraphics[width=0.33\linewidth]{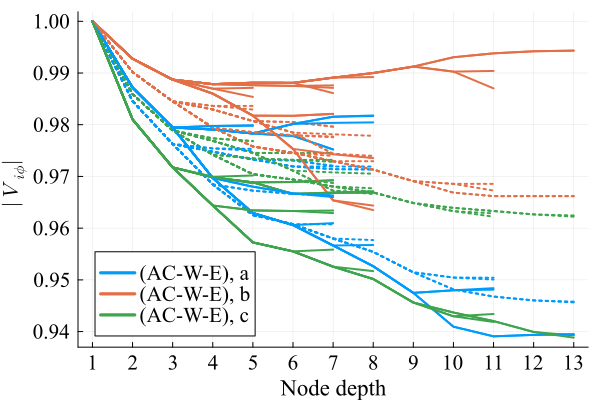}\label{fig:w-wye-delta}}
  \subfloat[][(AC-D)]{\includegraphics[width=0.33\linewidth]{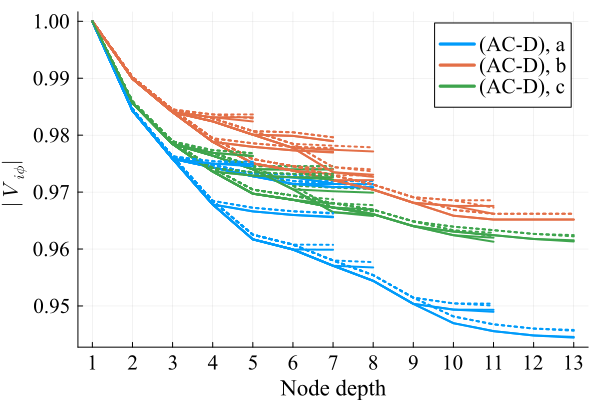}\label{fig:w-ac-ac-constp}}
  \subfloat[][(LP-D-E)]{\includegraphics[width=0.33\linewidth]{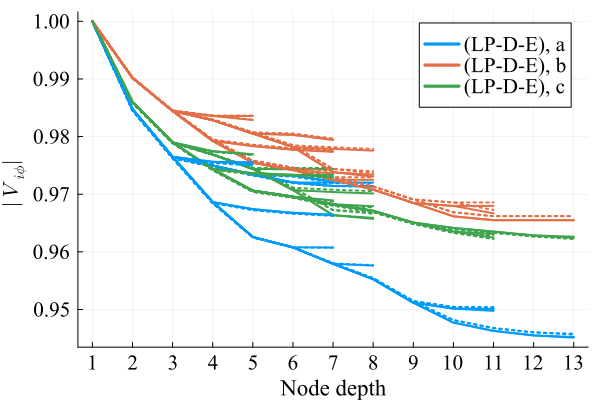}\label{fig:w-ac-lp}}
  \caption{{\color{black}Comparison of voltage magnitudes in p.u. for IEEE 37, with dotted lines representing (AC-D-E).}}
  \label{fig:wye-delta}
\end{figure*}
{\color{black}
We highlight the significance of integrating delta connections and voltage-dependent loads into OPF analyses. In this section, we focus on the IEEE 37-bus distribution feeder as it has only delta connected loads and half of them are voltage dependent. Consider two variants of the network: (AC-W-E) and (AC-D), where in (AC-W-E), delta-connected loads are treated as if they were connected in wye configurations, while in (AC-D), all loads are assumed to demand constant power. A comparison of the voltage magnitudes computed by these models with those of the base case (AC-D-E) reveals significant differences, as illustrated in Figures \ref{fig:w-wye-delta} and \ref{fig:w-ac-ac-constp}, respectively. Specifically, (AC-W-E) underestimates the voltage magnitude of phase b while overestimating that of phase a, whereas (AC-D) consistently underestimates voltage magnitudes across the network. However, when considering delta connections and voltage-dependent loads, even the linear approximation of the base case (AC-D-E) represented by (LP-D-E) yields a much-aligned voltage profile as (AC-D-E), as depicted in Figure \ref{fig:w-ac-lp}. It underscores the importance of considering delta connections and the voltage dependency of loads in OPF analyses.
}

{\color{black}
\subsection{Ramifications and Limitations of Assumption \ref{assum:linearized}}
The proposed linear model relies on Assumption \ref{assum:linearized}, which may not accurately represent unbalanced systems with line losses. In this section, we examine the implications of Assumption \ref{assum:linearized}. 
\begin{figure}[t!]
  \centering
  \includegraphics[width=0.6\linewidth]{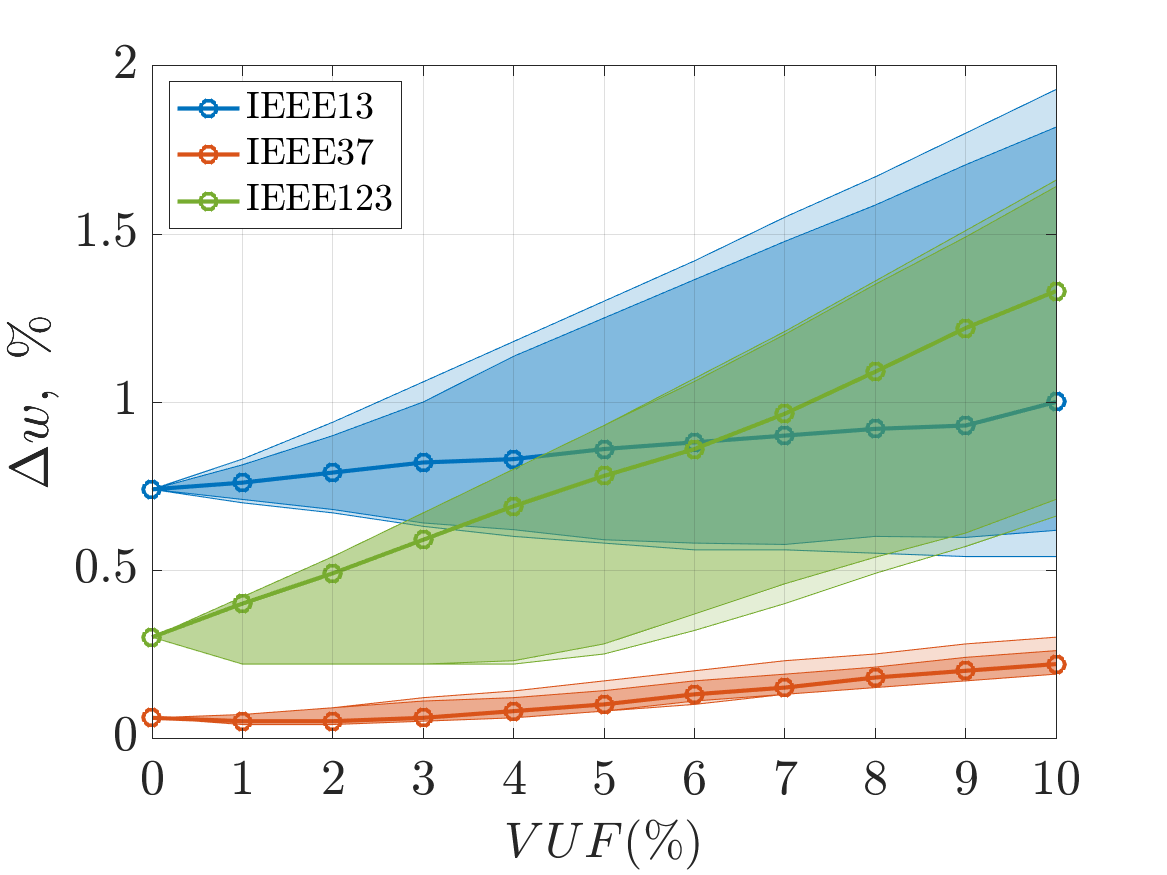}
  \caption{{\color{black}Impact of voltage imbalance.}}
  \label{fig:werror-imbal}
\end{figure}
To evaluate the impact of angle imbalance, we introduce perturbations to the phase angles of the substation reference voltage $V^{ref}_0$. We adhere to a formal definition of phase imbalance of a three-phase voltage $V=[V_a; V_b; V_c]$, referred to as $VUF$ \cite{von2001assessment}, defined as:
\begin{equation}
VUF(\%) = \frac{|V_n|}{|V_p|} \times 100,
\label{eq:vuf}
\end{equation}
where $V_p=\frac{V_a + \gamma^2 V_b + \gamma V_c}{3}$ and $V_n=\frac{V_a + \gamma V_b + \gamma^2 V_c}{3}$, with $\gamma = e^{-\mathrm{i} (2\pi/3)}$. For each $VUF$ value ranging from 1\% to 10\%, incremented by one, we randomly select a feasible angle perturbation on phase $c$ and compute the corresponding perturbation on the angle of phase $b$ so that the resultant $V_0^{ref}$ has the chosen $VUF$ value. This process is repeated 100 times for each $VUF$ value.

Figure \ref{fig:werror-imbal} illustrates the average relative difference in the squared voltage magnitude between (LP-D-E) and (AC-D-E) with respect to the varying levels of $VUF$; the shaded region represents the minimum, 10th, 90th percentiles, and maximum values, while the solid line with markers depicts the median values. The figure demonstrates a trend where the error increases as the level of imbalance grows. Nevertheless, it is noteworthy that for $VUF$ within 2\%, the error remains below 1\% for all test systems, and occurrences of larger $VUF$ values may be infrequent per several standards establishing permissible limits on $VUF$. For example, IEEE Standard 1159 \cite{475495} proposes a maximum limit of 2\% for both low- and medium-voltage networks, while similarly, IEC Standard 61000-3-13 advocates for a 2\% limit for low- and medium-voltage systems.


\begin{figure}[t!]
  \centering
  \includegraphics[width=0.6\linewidth]{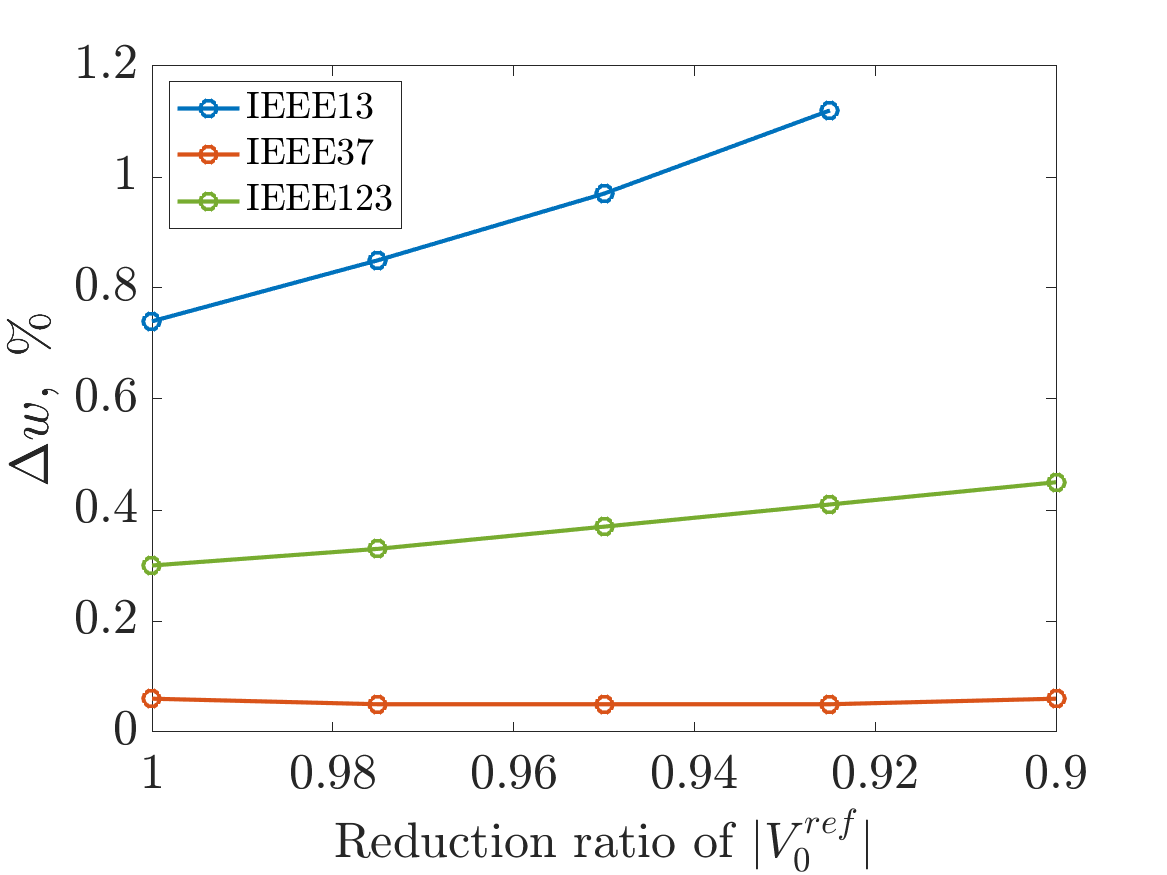}
  \caption{{\color{black}Impact of voltage magnitude reduction.}}
  \label{fig:werror-re}
\end{figure}
Subsequently, we investigate the effects of reducing the magnitude of $V^{ref}_0$ to assess the implications of ignoring line losses. We decrease the magnitude of the reference voltage by a factor $m$, ranging from 1 to 0.9 with decrements of 0.25. Figure \ref{fig:werror-re} displays the average relative difference in the squared voltage magnitude. As the magnitude decreases, line losses increase, leading to larger errors. Notably, IEEE 37 is less affected within the tested range, whereas IEEE 13 is more vulnerable to magnitude reduction, with the maximum error reaching around 1.2 for a reduction by a factor of 0.925. For the IEEE 13 system, a reduction ratio of 0.9 renders both (AC-D-E) and (LP-D-E) infeasible.
}
\section{Conclusions} \label{sec:conclusion}
We proposed a linear approximation of OPF for multiphase radial networks with mixed wye and delta connections as well as exponential load models. We proposed a system of linear equations that exactly illustrates the bus power withdrawal/injection based on the power flow from delta-connected devices under the assumptions made for the linear unbalanced OPF model proposed by \cite{gan2014convex}. The proposed system of linear equations can be used for various delta-connected devices, such as generators, loads, and capacitors. Numerical studies on IEEE 13-, 37-, and 123-bus systems showed that the linear approximation produces solutions with good empirical error bounds in a short amount of time, suggesting its potential applicability to various planning and operations problems with advanced features (e.g., line switching, capacitor switching, transformer controls, and reactive power controls). {\color{black}Furthermore, an experiment involving varied degrees of voltage angle imbalance and magnitude reduction underscored the effectiveness of our proposed model in maintaining acceptable error bounds, even within systems subject to practical deviations from the assumptions of small loss and balanced voltage. Future research endeavors will focus on improving performance under high voltage imbalances, possibly through the integration of data-driven line-specific $\Gamma$ matrices, alongside the inclusion of accurate modeling of three-phase transformers \cite{low2023three}, and comparison with successive approximation approaches \cite{bernstein2017linear,karagiannopoulos2018centralised}.}


\bibliographystyle{IEEEtran}
\bibliography{References.bib}

\end{document}